\documentclass[12pt,a4paper,oneside,english]{amsart}
\usepackage{graphicx}
\usepackage{setspace}
\usepackage{amssymb}
\makeatletter
\usepackage{babel}

 \theoremstyle{plain}
 \newtheorem{thm}{Theorem}[section]
 \numberwithin{equation}{section} 
 \numberwithin{figure}{section} 
 \theoremstyle{plain}
 \newtheorem*{thm*}{Theorem}
 \theoremstyle{definition}
 \newtheorem*{defn*}{Definition}
 \theoremstyle{plain}
 \newtheorem{prop}[thm]{Proposition} 
 \theoremstyle{remark}
 \newtheorem{rem}[thm]{Remark}
 \theoremstyle{remark}
 
\theoremstyle{plain}
 
 \theoremstyle{plain}
 \newtheorem{lem}[thm]{Lemma} 
 \theoremstyle{definition}
 \newtheorem{defn}[thm]{Definition}
\newtheorem*{acknowledgment}{Acknowledgment}

\AtBeginDocument{
  
}

\begin{document}

\title{Rewriting systems and embedding of monoids in groups}
\author{Fabienne Chouraqui}

\begin{abstract}
 In this paper, a  connection between rewriting systems and
embedding of monoids in groups is found. We  show that if a group with a positive presentation has  a  complete rewriting system $\Re$ that satisfies the
condition that each rule in $\Re$ with positive left-hand side has
a positive right-hand side, then the monoid presented by the subset of positive rules from $\Re$  embeds in the group. As an example, we give a simple proof that right angled Artin
 monoids embed in the corresponding right angled  Artin
 groups. This is a special case of the well-known result of Paris \cite{paris}
  that Artin monoids embed in their groups.
\end{abstract}
\maketitle
\section{Introduction}
In \cite{adian}, Adian gave a sufficient condition on the presentation of a monoid
 for its  embeddability  in a group, which is very simple to check.

 Adian's results have been extended in different directions in the
 works
of several authors, for example, in the paper of Kashintsev
\cite{katshin} where various small cancellation conditions are
used to study embeddability of semigroups in groups and  in the
paper of Krstic \cite{krstic}.

 In this work, instead
  looking  at  the presentation of the monoid, we consider the group presentation.
 In section $2$, we give some preliminaries on  embeddability  of monoids in groups
 and on rewriting systems.
In section $3$, we give the main result and we give some
definitions.
 We define  a positive word in the
  group $G=\operatorname{Gp}\langle X\mid R\rangle  $  to be a word which  belongs
to the free monoid $X^{*}$ and we consider the empty word $``1''$ as
 a positive word. We say that a  rewriting system $\Re$ for the group
  $G=\operatorname{Gp}\langle X\mid R\rangle  $ satisfies the condition $C^{+}$
  if each rule in $\Re$
 with positive left-hand side has
 a positive right-hand side. We define $\Re^{+}$ to be the subset of all rules of $\Re$
with positive left-hand side. If $\Re$  satisfies the condition
$C^{+}$ then the rules in $\Re^{+}$  are called positive rules.
\pagebreak

The main result proved in section $3$ can then be stated as follows:
\begin{thm*}
Let $N$ be a monoid  and $G = \operatorname{Gp}\langle X\mid R\rangle  $ a group. Assume that $G$  has
a  complete rewriting system $\Re$ that satisfies the condition
$C^{+}$.  Then the monoid $M=\operatorname{Mon}\langle X\mid \Re^{+}\rangle  $ embeds into $G$ as the monoid of positive words relative to this presentation. In particular, if $N$  is isomorphic to $M$, then  $N$ embeds into $G$.
\end{thm*}

 Using this result, in section $4$, we give a simple and natural proof  that right angled Artin monoids
  embed in their corresponding right angled  Artin
 groups. This is a special case of the well-known
 result of Paris \cite{paris}
 that Artin monoids embed in their groups. If $G$ is a right angled  Artin
 group, then we define for $G$ a terminating rewriting system $\Re_{0}$,
 that is based on a length-lexicographical ordering induced by a  total ordering on the generators,
  which satisfies the condition $C^{+}$ and such that $\Re_{0}^{+}$ is a rewriting system for the monoid with the same presentation.
 We then apply the Knuth-Bendix algorithm of completion on $\Re_{0}$ in order
 to obtain a complete (maybe infinite) rewriting system $\Re$, which is equivalent to $\Re_{0}$.
 We then show that  the rewriting system $\Re$ satisfies the condition
 $C^{+}$ and we show that at each step of the completion a positive rule is  obtained from
 an ambiguity between positive rules which were obtained at a
 former step.\\
 This gives an algorithm for checking the embeddability of a
 monoid $M$ generated by a set $X$ into a group $G = \operatorname{Gp}\langle X\mid R\rangle  $ which can be described in the
 following way:
Define for $G$ a terminating rewriting system $\Re_{0}$,
 which is based on a total ordering of the words,  that satisfies the condition
 $C^{+}$ and such that $\Re_{0}^{+}$ is a rewriting system for $M$.\\
 Apply the Knuth-Bendix algorithm of completion on $\Re_{0}$ in order
 to obtain a complete  rewriting system $\Re$, which is equivalent to
 $\Re_{0}$.\\
 If the complete  rewriting system $\Re$ does not satisfy the condition
 $C^{+}$, then nothing can be concluded.\\
If the complete  rewriting system $\Re$ does  satisfy the
condition  $C^{+}$ and at each step a  positive rule is  created by an ambiguity
between positive rules which were created at a  former step, then $M$ embeds in $G$.

\begin{acknowledgment}
This work is a part of   my Ph-d thesis, done  at the Technion under
the supervision of Professor Arye Juhasz.  I am very grateful to
Professor Arye Juhasz, for his patience, his encouragement and his many judicious remarks and  to the referee for his very useful comments.
\end{acknowledgment}

\pagebreak
\section{Preliminaries}
\subsection{On embedding of monoids in groups}\
Let $\Sigma$ be a non-empty set. We denote by $\Sigma^{*}$the free
monoid generated by $\Sigma$; elements of $\Sigma^{*}$ are finite
sequences called \textbf{words} and the empty word will be denoted
by 1.
Let $M$ be a monoid. The left and right  laws of cancellation hold in $M$ if the
following conditions are satisfied: if $ab=ac$ implies $b=c$ and
 $ba=ca$ implies $b=c$ respectively. A commutative monoid with laws of
cancellation can always be embedded in a group. This is not
necessarily true for non-commutative monoids, so for this class of
monoids the question if there is a group $G$ such that $M$ can be
embedded in $G$ arises. Adian, in \cite{adian}, gives a partial answer to
this question: he gives a sufficient condition for embeddability
of the monoid in the group with the same presentation. In what
follows, we will describe Adian's criteria for embeddability.
Let $M$ be a monoid with presentation $\langle \Sigma\mid
r_{1},r_{2},..,r_{m}\rangle  $ , where each relation $r_{i}$ have the form
$u_{i}=v_{i}$, $u_{i},v_{i} \in \Sigma^{*}$ for $1\leq i\leq m$ ,
and none of the words $u_{i},v_{i}$ is empty. Adian in
\cite{adian} defined for such a presentation the following two
graphs, which are called the left and right graphs of the
presentation. In the \textbf{left graph, $\Gamma_{l}$}, to each
relation $r_{i}$ there corresponds an edge
\textbf{$(a_{i}^{1},b_{i}^{1})$} , where \textbf{$a_{i}^{1}$} is
the first letter of the word $u_{i}$ and \textbf{$b_{i}^{1}$} is
the first letter of the word $v_{i}$. In the \textbf{right graph,
$\Gamma_{r}$}, to each relation $r_{i}$ there corresponds an edge
\textbf{$(a_{i},b_{i})$} , where \textbf{$a_{i}$} is the last
letter of the word $u_{i}$ and \textbf{$b_{i}$} is the last letter
of the word $v_{i}$. A sequence of edges
$(a_{1},a_{2}),(a_{2},a_{3}),..,(a_{k-1},a_{k})$ is called a path
connecting $a_{1}$ and $a_{k}$ and if $a_{1}=$$a_{k}$ then this
path is closed. A closed path of length greater than 1 such that
the first and last edges of the path are distinct, i.e. correspond
to different relations and all vertices $a_{j}$, for $1<j<k$, are
different from $a_{1}$ is called a cycle. A closed path of length
1 is an elementary cycle.
\begin{thm}
\cite[p.57]{adian} Let $M$ be a monoid with presentation\\
$\langle \Sigma\mid r_{1},r_{2},..,r_{m}\rangle  $ , where each relation $r_{i}$
have the form \ $u_{i}=v_{i}$,  $u_{i},v_{i} \in \Sigma^{*}$ for
$1\leq i\leq m$,
 and none of
the words $u_{i},v_{i}$ is empty. Assume there are no cycles in
the left and right graphs of the presentation of $M$ and let
$w_{1},w_{2}$ be words in $M$. Then the equality $w_{1}=w_{2}$
holds in the monoid $M$ if and only if it holds in the group $G$
presented in terms of the same generators and the same defining
relations and in this case $M$  embeds in $G$.
\end{thm}
\textbf{Examples}
\begin{itemize}
\item[1.] Let $M=\langle a,b\mid ab=ba\rangle  $.The left and right graphs $\Gamma_{l}$
and $\Gamma_{r}$ are the same and none of them has any cycle, so
the monoid embeds in the free abelian group on two
generators.
\item[2.]
Let $M=\langle x,y\mid xy=y\rangle  $. There is no left cycle but there is an
elementary right cycle (from $y$ to itself). That means that the monoid is
left cancellative but not necessarily right cancellative.
\end{itemize}

\subsection{On rewriting systems}\

 We refer the
reader to \cite{book}, \cite{cohen} and \cite{huet} for more
details.
\begin{defn*}
A \textbf{rewriting system} $\Re$ on $\Sigma$ is a set of ordered
pairs in $\Sigma^{*}\times \Sigma^{*}$ .
\end{defn*}
If $(l,r) \in \Re$ then for any words $u$ and $v$ in $\Sigma^{*}$,
 we say the word $ulv$ \textbf{reduces} to the word $urv$ and we
write $ulv\rightarrow urv$ . A word $w$ is said to be \textbf{reducible}
if there is a word $z$ such that $w\rightarrow z$. If there is no
such $z$ we call $w$ \textbf{irreducible}.

A rewriting system $\Re$ is called \textbf{terminating} \textbf{(or
Noetherian)} if there is no infinite sequence of reductions
$w_{1}\rightarrow w_{2}\rightarrow...\rightarrow w_{n}\rightarrow....$

We denote by  ``$\rightarrow^{*}$''  the reflexive
transitive closure of the relation $``\rightarrow''$.

A rewriting system $\Re$ is called \textbf{Church - Rosser} if for
any words $u,v$ in $\Sigma^{*}$, $u\leftrightarrow^{*}v$ implies
that there is a word $z$ in $\Sigma^{*}$ such that
$u \rightarrow^{*}z$ and $v \rightarrow^{*}z$ (i.e. if $u$ and $v$
are equivalent then they have a common descendant). $\Re$ is
called \textbf{confluent} if for any words $u,v,w$ in $\Sigma^{*}$ ,
$w\rightarrow^{*}u$ and $w\rightarrow^{*}v$ implies that there is
a word $z$ in $\Sigma^{*}$ such that $u\rightarrow^{*}z$ and
$v\rightarrow^{*}z$ (i.e. if $u$ and $v$ have a common ancestor
then they have a common descendant). These two properties of $\Re$
are equivalent: \textbf{$\Re$ is Church-Rosser if and only if
$\Re$ is confluent}\cite{book}. For any word $w$ in $\Sigma^{*}$,
the confluence of $\Re$ ensures the existence of at most one
irreducible equivalent word.

$\Re$ is called \textbf{locally confluent} if for any words $u,v,w$
in $\Sigma^{*}$ , $w\rightarrow u$ and $w\rightarrow v$ implies that
there is a word $z$ in $\Sigma^{*}$ such that $u \rightarrow^{*} z$
and $v \rightarrow^{*} z$ .

$\Re$ is called \textbf{complete (or convergent)} if $\Re$ is
terminating and confluent. So, if $\Re$ is complete then every
word $w$ in $\Sigma^{*}$ has a unique irreducible equivalent word
$z$ and  $z$ is called the \textbf{normal form} of $w$.

It is a very hard task to determine the completeness of an arbitrary
rewriting system. Knuth and Bendix have elaborated an algorithm which
for a given finite and terminating rewriting system $\Re$, tests
its completeness and if $\Re$ is not complete then new rules are
added to complete it. Instead of testing the confluence of $\Re$,
the algorithm tests the locally confluence of $\Re$, since \textbf{for
a terminating rewriting system $\Re$ locally confluence and confluence
are equivalent} \cite[p.16]{book}.

We say that two rewriting systems $\Re$ and $\Re'$ are
\textbf{equivalent} if :\
$w_{1}\leftrightarrow^{*}w_{2}$ modulo
$\Re$ if and only if $w_{1}\leftrightarrow^{*}w_{2}$ modulo
$\Re'$. So, by applying the Knuth-Bendix algorithm on a
terminating rewriting system $\Re$ a complete rewriting system
$\Re'$ which is equivalent to $\Re$ can be found. If the algorithm
halts after a finite number of steps and the rewriting system
obtained $\Re'$ is complete  then $\Re'$ is finite and complete,
since a finite number of rules are added.

We call the triple of non-empty words $u,v,w$ in $\Sigma^{*}$ an \textbf{overlap
ambiguity} if there are $r_{1},r_{2}$ in $\Sigma^{*}$ such that $uv\rightarrow r_{1}$
and $vw\rightarrow r_{2}ý$ are rules in $\Re$. We then say that $r_{1}w$
and $ur_{2}$ are the corresponding \textbf{critical pair.} When the triple $u,v,w$
 in $\Sigma^{*}$ is an overlap ambiguity
we will say that the rules $uv\rightarrow r_{1}$ and $vw\rightarrow r_{2}$
\textbf{overlap at $v$} or that there is an \textbf{overlap}
between the rules $uv\rightarrow r_{1}$ and $vw\rightarrow r_{2}$.
If there exists a word $z$ such that $r_{1}w\rightarrow^{*}z$ and
$ur_{2}\rightarrow^{*}z$ , then we say that the critical pair resulting
from the overlap of the rules \textbf{}$uv\rightarrow r_{1}$
and $vw\rightarrow r_{2}$ in $\Re$ \textbf{resolves}.

The triple $u,v,w$ of possibly empty words in $\Sigma^{*}$ is called
an \textbf{inclusion ambiguity} if there are $r_{1},r_{2}$ in
$\Sigma^{*}$(which must be distinct if both $u$ and $w$ are empty,
but otherwise may be equal) such that $v\rightarrow r_{1}$ and
$uvw \rightarrow r_{2}$ are rules in $\Re$. We then say that
$ur_{1}w$ and $r_{2}$ are the corresponding \textbf{critical
pair}. If there exists a word $z$ such that
$ur_{1}w\rightarrow^{*}z$ and $r_{2}\rightarrow^{*}z$ , then we
say that the critical pair resulting from the inclusion ambiguity
of the rule $v\rightarrow r_{1}$ in $uvw\rightarrow r_{2}$ in
$\Re$ \textbf{resolves}. Given a terminating rewriting system
$\Re$, the Knuth-Bendix algorithm may halt after a finite number
of steps and fail if the words in a critical pair cannot be
compared in order to resolve it. If a total ordering of the words
is defined, then  the algorithm ensures the existence of a
complete rewriting system  $\Re'$ which is equivalent to $\Re$
(see \cite{huet}).

 So, $\Re$ is  \textbf{complete} if
$\Re$ is terminating and locally confluent or in other words if
$\Re$ is terminating and all the critical pairs resolve.

\section{The main result: a criteria for embeddability of monoids in groups}

In this section, we use the following notation.\\
Let $M$ be a monoid  and let $G$ be the group presented by $\operatorname{Gp}\langle X\mid R\rangle  $, so the monoid presentation of $G$ is
$\operatorname{Mon}\langle X\cup X^{-1} \mid R\cup R_{0}\rangle  $ , where $X^{-1}$ denotes the
set $\{x_{0}^{-1},...,x_{n}^{-1}\}$ and $R_{0}=\{ x_{i}x_{i}^{-1}=1,x_{i}^{-1}x_{i}=1,$ with
 $i\in\{0,1,...,n\}\}$.
We refer the reader to \cite{sims} for more details.
\begin{defn}
Let $w$ be a word in $(X\cup X^{-1})^{*}$. We say that $w$ is a
\textbf{positive word in $G$} if $w$ belongs to the free monoid
$X^{*}$. The empty word $``1''$ is considered here as a positive
word.
\end{defn}
\textbf{Example: The braid group $B_{3}$}
 with the presentation $B_{3}= \operatorname{Gp}\langle a,b,c \mid a^{3}=b^{2}=c \rangle  $.\\
The word $ac$ is a positive word in $B_{3}$, while the word $a^{-1}c$ is not.
\begin{defn}
We denote by \textbf{$\Re^{+}$} the subset of all rules of $\Re$
with positive left-hand side. We allow  in $\Re^{+}$  rules of the
form $l \rightarrow 1$, where $l$ is a positive word.
\end{defn}
\begin{defn}
Let  $\Re$ be a rewriting system for $G$, with $\Re^{+}\neq \emptyset$.
We say that \textbf{$\Re$ satisfies the condition $C^{+}$} if
each rule in $\Re$ with positive left-hand side has  a positive
right-hand side.
\end{defn}

\textbf{Example:} In \cite{pederson}, the authors give a finite and
complete rewriting
system for the group $B_{3}$, which is described as follows:\\
$\Re = \{a^{-1}\rightarrow c^{-1}a^{2} ,
 b^{-1}\rightarrow c^{-1}b , a^{3} \rightarrow c,
b^{2} \rightarrow c , ac\rightarrow ca ,
ac^{-1}\rightarrow c^{-1}a , bc\rightarrow cb ,
bc^{-1}\rightarrow c^{-1}b , cc^{-1}\rightarrow 1 ,
c^{-1}c\rightarrow 1 \}$\\
The rewriting system $\Re$  satisfies the condition  $C^{+}$, since  each rule in
$\Re$ with positive left-hand side has  a positive
right-hand side.
The rewriting system $\Re^{+}$ is: $\Re^{+} = \{ a^{3} \rightarrow c,
b^{2} \rightarrow c , ac\rightarrow ca ,
bc\rightarrow cb \}$
\begin{rem}
One can check that the rewriting system $\Re^{+}$ described above is a finite and
 complete rewriting system
for the monoid presented by $\langle a,b,c \mid a^{3}=b^{2}=c \rangle  $.
\end{rem}
\begin{defn}
We say that $\Re^{+}$ is a rewriting system for $M$ if $M$ is
isomorphic to the factor monoid $X^{*}/\equiv_{\Re^{+}}$.
\end{defn}
\begin{lem}
\label{lem:R+_complete}  Assume that $G$ has a complete rewriting system $\Re$ which
satisfies the condition $C^{+}$. If $\Re^{+}$ is a rewriting system for $M$, then $\Re^{+}$ is a  complete
rewriting system for $M$. Further, if $\Re$  is finite then $\Re^{+}$ is
also finite.
\end{lem}

\begin{proof}
We have to show that  $\Re^{+}$ is terminating and that all the critical pairs
resulting from any kind of ambiguity between rules in $\Re^{+}$ resolve by using
rules in $\Re^{+}$.\\
Assume there is an infinite sequence of reductions in $\Re^{+}$: \\
$w_{1} \rightarrow w_{2} \rightarrow ... w_{m} \rightarrow ..\rightarrow ..$\\
Since the rules applied in this sequence are rules in  $\Re^{+}$,
these are also rules in  $\Re$. So, there is an infinite sequence
of reductions in $\Re$. But this is a contradiction to the fact
that $\Re$ is terminating. So, $\Re^{+}$ is terminating.\

Assume that $u$ and $v$ is a critical pair resulting from any kind of ambiguity
 between rules in $\Re^{+}$.
 Since $\Re$ is complete, this critical pair resolves using rules in $\Re$.
But the words $u$ and $v$ are positive words, since $\Re$ satisfies the condition $C^{+}$.
So, the rules from $\Re$  used to resolve this critical pair belong to $\Re^{+}$.
So, $\Re^{+}$ is complete.
If  $\Re$ is finite, then  $\Re$ is also finite, since  $\Re^{+}$ is a subset of  $\Re$.
\end{proof}

\begin{defn}
Let $\Re$ be a rewriting system. We denote by $\equiv_{\Re}$ the
 equivalence relation generated by $\Re$. In fact, $\equiv_{\Re}$ is a congruence since it is compatible with concatenation. (see \cite{sims})
\end{defn}

\begin{lem}
\label{cor:same_irred}  Assume that $G$ has a complete rewriting system $\Re$ which
satisfies the condition $C^{+}$. Let $u,v$ be positive words. Then
$u\equiv_{\Re} v$ if and only if $u\equiv_{\Re^{+}} v$.
\end{lem}

\begin{proof}
If $u\equiv_{\Re^{+}} v$, then clearly $u\equiv_{\Re} v$, since
$\Re^{+} \subseteq \Re$.\\
Assume that  $u\equiv_{\Re} v$. Then there is a (modulo $\Re$) irreducible word
 $z$ such that $u \rightarrow^{*} z$ and $v
\rightarrow^{*} z$. The words $u$ and $v$ are positive words in
$G$ and $\Re$ satisfies the condition $C^{+}$, so the rules used
in order to reduce $u$ and $v$ to $z$ belong to $\Re^{+}$ and  $z$
is also a positive word. So, $u\equiv_{\Re^{+}} v$.
\end{proof}

\begin{thm}\label{thm_result}
Let $N $ be a monoid  and  $G = \operatorname{Gp}\langle X\mid R\rangle  $ a group. Assume that $G$ has
a  complete rewriting system $\Re$ that satisfies the condition
$C^{+}$.  Then the monoid $M=\operatorname{Mon}\langle X\mid \Re^{+}\rangle  $ embeds into $G$ as the monoid of positive words relative to this presentation. In particular, if $N$  is isomorphic to $M$, then  $N$ embeds into $G$.
\end{thm}
\begin{proof}
Assume that $\Re^{+}$ is a rewriting system for $M$. Then from
lemma \ref{lem:R+_complete}, $\Re^{+}$ is a complete rewriting
system for $M$.
Let $u$ and $v$ be positive words such that $u$ and $v$ are equal in $G$.
 We  show that the words $u$ and $v$ are equal in $M$.
Since $\Re$ is a complete rewriting system for $G$, $u$ and $v$
are equal in $G$ implies that $u\equiv_{\Re} v$. From
lemma \ref{cor:same_irred}, $u\equiv_{\Re^{+}} v$ since $u$ and $v$
are positive words. But $\Re^{+}$ is a complete rewriting system
for $M$, so $u$ and $v$ are equal in $M$.\\
\end{proof}

\section{An Application:Embedding of right angled  Artin monoids}

 Recall that a group G
is called an Artin group and the corresponding monoid is called an Artin monoid
if it is presented by a set $X=\{x_{1},...,x_{n}\}$  subject to relations of the form
 $\underbrace {x_{i}x_{j}x_{i}..}_{m_{i,j}}=\underbrace{ x_{j}x_{i}x_{j}..} _{m_{j,i}}$
 where  $m_{i,j}= m_{j,i} \geq 2$ or $m_{i,j}=0 $ in which case we omit the relation between
$x_{i}$ and $x_{j}$.
 A right-angled Artin group (or partially commutative group) is one
  in which $m_{i,j} \in \{2,0 \}$. In other words, in the presentation for the Artin group, all relations are commutator
relations: $x_{i}x_{j}=x_{j}x_{i}$.
A standard way to specify the presentation for any Artin group
is by means of the defining graph $\Im$. Let $\Im$ be a finite, simplicial graph with edges
 labeled by integers greater than one. The Artin group (resp. the Artin monoid) associated to $\Im$ is the group (resp. monoid)  whose presentation has
generators  corresponding to the vertices $x_{1},.., x_{n}$ and there is a relation
$\underbrace {x_{i}x_{j}x_{i}..}_{m_{i,j}}=\underbrace{ x_{j}x_{i}x_{j}..} _{m_{i,j}}$
  for each edge labeled $m_{i,j}$ connecting  $x_{i}$ and $x_{j}$.
 Several authors have studied the question of embeddability of general Artin
monoids in Artin groups: for example special cases of this question have been
considered by Charney \cite{charney} and Cho and Pride \cite{cho}. Paris \cite{paris}
 proved the following  general result about Artin
groups and Artin monoids:  Every Artin monoid embeds in the corresponding Artin group.

 In this section, we show the following theorem.
 \begin{thm} \label{thm:artin_embed}
 Right angled  Artin monoids  embed in their corresponding groups.
 \end{thm}
 The proof is by showing that they admit a
   complete rewriting system $\Re$  satisfying the condition $C^{+}$
   and the condition on  $\Re^{+}$.

In \cite{hermiller2} Hermiller and Meier  constructed finite and
complete rewriting systems for graph groups, and as a special case, for right angled Artin groups using  another
presentation. In \cite{vanwyk}, Van Wyk constructs complete rewriting systems
for graph groups. His main interest is on the  normal forms.
  Although we obtain the same normal forms as in \cite{vanwyk},
  our main interest is on the rewriting system and on how the rules look like in the light of theorem \ref{thm_result}.

Let $\Im$ be the defining graph of a right-angled Artin group $G$.
Denote by $X$ the set of vertices in $\Im$ (or equivalently the
set of generators of $G$) and by $X^{-1}$ the set of inverses of
the elements in $X$.

\begin{lem}\label{lem:sym{ab=ba}}
To each edge $a - b$ in $\Im$,  there corresponds the following
set of $4$ relations in the monoid presentation of the group:\\
$\{a^{\epsilon}b^{\delta} = b^{\delta}a^{\epsilon}$, where $\epsilon,\delta$ take all values $\pm 1 \}$
\end{lem}
\begin{proof}
To each edge $a - b$ in $\Im$,  there corresponds the relation
$ab=ba$ in the group presentation. In the monoid presentation of
the group,  by adding the relations: $xx^{-1}=1$ and $x^{-1}x=1$
for every $x \in X$, we obtain a set with  many relations which
are all consequences of the set of $4$ relations described above.
\end{proof}

\begin{lem}\label{lem:order_resolve}
Each  set of $4$ relations corresponding to an edge $a - b$ in
$\Im$, as defined in lemma \ref{lem:sym{ab=ba}}, can be oriented
in two ways such that the overlaps with the rules  $xx^{-1}
\rightarrow 1$ and $x^{-1}x \rightarrow1$ $(x \in X)$ resolve:
$\{a^{\epsilon}b^{\delta} \rightarrow b^{\delta}a^{\epsilon}$, where $\epsilon,\delta$ take all values $\pm 1 \}$  or  $\{b^{\delta}a^{\epsilon} \rightarrow a^{\epsilon}b^{\delta}$, where $\epsilon,\delta$ take all values $\pm 1 \}$
\end{lem}
\begin{rem}\label{rem:R0satisfairC+}
Note that in both cases, when the left-hand side of a rule  is a
positive word its right-hand side is also a positive word.
\end{rem}
\begin{proof}
We will check that the overlap of the first rule (in the first set) with a rule of free
 reduction resolve and we omit this  technical checking for the other rules.
$$\begin{array}{ccccccccc}
 &&   \overbrace{a^{-1}a} b=a^{-1} \overbrace{ab} \\
 & \swarrow    & \searrow \\
b &&& \underbrace{a^{-1}b} a \rightarrow b \underbrace{a^{-1}a} \rightarrow b\\
\end{array}$$
\end{proof}

From lemma \ref{lem:order_resolve}, we have that for any right angled Artin group $G$, by defining a total ordering on the generators and using the  length-lexicographical ordering induced by it, one can define a rewriting system of the following form: $\Re_{0} =\{ x^{\gamma}y^{\epsilon} \rightarrow y^{\epsilon}x^{\gamma}$,  where $\gamma$ and $\epsilon$ take all values $\pm 1$, for each edge $ x-y \in \Im \}$.
If  $\Re_{0}$ is finite and complete with no need of completion, then it satisfies the condition $C^{+}$ (see remark \ref{rem:R0satisfairC+}) and $\Re_{0}^{+}$ is a rewriting system for $M$. Yet, this is not true in
general, an example of infinite complete rewriting system is given
in \cite[p.250]{hermiller2}. We will show here (propositions \ref{prop:Rsatisfait_pos} and \ref{prop:RsatisfiesR+}) that by applying the Knuth-Bendix algorithm of completion on $\Re_{0}$, an
equivalent rewriting system (not necessarily finite) $\Re$ is
obtained which is complete and which satisfies  the condition
$C^{+}$ and the condition on $\Re^{+}$. Since $\Re_{0}$ is
terminating and there is a total ordering of the words, the
algorithm ensures the existence of an equivalent rewriting system
$\Re$ which is complete. This will imply that the
monoid embeds in the corresponding group, by theorem \ref{thm_result}.
\begin{rem}
 In order to avoid inclusion ambiguities,
we will use the following strategy of reduction in the
Knuth-Bendix algorithm of completion: at each step, we reduce each
word in the critical pair  to its normal
form modulo the  rewriting system obtained at the earlier step. So, we need consider only overlap ambiguities.
\end{rem}
Before we proceed, we need the following definitions and notations.\\
Let $w=t_{1}t_{2}...t_{k}$, where $t_{i} \in X \bigcup X^{-1}$ for $1 \leq i \leq k$.\\
We define the \textbf{prefix of $w$}
  to be the following set of words: \\
   $\operatorname{pref}(w)=\{ t_{1}, t_{1}t_{2}, t_{1}t_{2}t_{3} ,.., t_{1}t_{2}t_{3}...t_{k} \}$ \\
   We denote by $\ell(w)$ the length of the word $w$. \\
   We denote by $\Re$ the complete rewriting system obtained by applying the
   Knuth-Bendix algorithm of completion on $\Re_{0}$
    and we denote by $\Re_{n}$ the  rewriting system  obtained at the
$n-$th step in the algorithm of completion,  that is
$\Re_{n}=\Re_{n-1} \bigcup \{$the rules obtained at the $n-$th step$\}$. \\
We use the following definition from
\cite{vanwyk}: the set $\{t_{1},t_{2},..,t_{k}\}$ is defined to be
a \textbf{clique} if for all $1\leq i,j \leq k$, $i \neq j$, $t_{i}$
and $t_{j}$ commute.
\begin{lem}\label{lem:petitrule impliq}
Assume the  following set  of rules occur in $\Re_{0}$:
 $$(*)\left\{\begin{array}{c}
t_{1}x \rightarrow xt_{1}  \\
xt_{i} \rightarrow t_{i}x \\
\end{array}\right.$$
where    $t_{1},t_{i},x \in X \bigcup X^{-1}$ for $2 \leq i \leq k$ and such that
 $\{t_{1},t_{2},..,t_{k},x\}$ is not a clique.\\
 Assume there
is a rule $ux \rightarrow xu$ in $\Re$, where  $u
=t_{1}t_{i_{2}}..t_{i_{k}}$   and \\$t_{i_{j}} \in \{t_{2},..,t_{k}\}$. Then \\
(i) the
following  rules also occur in $\Re$: $u'x \rightarrow xu'$ for every $ u' \in
\operatorname{pref}(u)$. Further,  these
rules  are obtained at a step earlier than $ ux \rightarrow xu$.\\
(ii) the following rule also occurs in $\Re$: $ux^{-1} \rightarrow x^{-1}u$
\end{lem}
The proof is  technical and appears in appendix A. The main idea is that we apply step by step the algorithm of completion on the rules in (*) and the existence of the rule $ux \rightarrow xu$ in $\Re$ ensures the creation of the rules $u'x \rightarrow xu'$ for every $u' \in
\operatorname{pref}(u)$ and the rule $ux^{-1} \rightarrow x^{-1}u$.\\

For brevity of notation, we say that
\textbf{the  rule $ux \rightarrow xu$ satisfies the condition on prefixes}
if $u=t_{1}t_{2}..t_{k}$, where $x,t_{i} \in X \bigcup X^{-1}$
and the rules $t_{1}x \rightarrow xt_{1}$ and
  $xt_{i} \rightarrow t_{i}x$ for
 $2 \leq i \leq k$ occur in $\Re_{0}$.

\begin{defn}
Let $uv \rightarrow u'$ and $vw \rightarrow v'$ be rules in a
rewriting system  $\Re$, where $v$ is not the empty word. Then we
say that there is an overlap between  these rules with \textbf{$uv
\rightarrow u'$ at left and $vw \rightarrow v'$ at right}. We call
$v$ the \textbf{overlapping word} and $\ell(v)$ \textbf{the length
of the overlap}.
\end{defn}
In the following lemma, we  show that the length of the overlaps in
$\Re$ cannot exceed one, that is the overlapping words belong to
$X\bigcup X^{-1}$.
\begin{lem}
 Assume the rules   $ux_{1}x_{2}..x_{m} \rightarrow
 x_{m}ux_{1}x_{2}..x_{m-1}$ and $\allowbreak x_{1}x_{2}..x_{m}y \rightarrow
 yx_{1}x_{2}..x_{m}$   satisfy the condition on prefixes and occur in  $\Re$. Then $m=1$.
\end{lem}

\begin{proof}
Assume $m > 1$.
The rule $ux_{1}x_{2}..x_{m} \rightarrow
 x_{m}ux_{1}x_{2}..x_{m-1}$  satisfies the condition on  prefixes, so the  rule $x_{m}x_{1} \rightarrow x_{1}x_{m}$ occurs in $\Re_{0}$, that is (*) $x_{m}>x_{1}$.
 The rule $x_{1}x_{2}..x_{m}y \rightarrow
 yx_{1}x_{2}..x_{m}$   satisfies the condition on  prefixes, so  the  rule $yx_{m} \rightarrow x_{m}y$  occurs  in $\Re_{0}$, that is  $x_{1}>y$ and  $y>x_{m}$. So, we have $x_{1} > y > x_{m}$. But, this
 contradicts $(*)$, so $m=1$.
\end{proof}

\begin{prop} \label{prop:Rsatisfait_pos}
(i) The equivalent complete rewriting system $\Re$ obtained from the
Knuth-Bendix algorithm
 of completion of $\Re_{0}$
satisfies the following conditions:\\
 (A) All the rules in $\Re$ have the following form $ux \rightarrow
 xu$, \\where $u=t_{1}t_{2}..t_{k-1}t_{k}$ is a word, $x, t_{i} \in X\bigcup X^{-1}$, and they satisfy the condition on prefixes. \\
  (B) Only overlaps between rules from $\Re$ at left and rules
  from $\Re_{0}$ at right produce new rules, that is all the other
  kinds of overlaps resolve.\\
  (C)At the $n-$th  step of completion, the rules produced have
   words of length $n+2$ in each side.\\
(ii) $\Re$ satisfies the following condition: one side of
a rule in $\Re$ is a positive word  if and only if the other side is a positive word.
 In particular,  $\Re$ satisfies the  condition  $C^{+}$ (even if it is not finite).
\end{prop}

\begin{proof}
(i) The proof will be by induction on the number of steps in the
Knuth-Bendix algorithm of completion of $\Re_{0}$.
We show by induction that $\Re$ satisfies the conditions (A), (B) and (C).\\

At the $0-$th step, for  each rule of kind $x^{\gamma}y^{\epsilon}
\rightarrow  y^{\epsilon}x^{\gamma}$ in $\Re_{0}$, where
$\gamma,\epsilon = \pm1$,  assumption  (A) holds trivially. Each
rule has words of length $2$ in each side.\\
At the first step, we have overlaps between  rules
 $x^{\gamma}y^{\epsilon} \rightarrow  y^{\epsilon}x^{\gamma}$ and
 $y^{\epsilon}z^{\delta} \rightarrow  z^{\delta}y^{\epsilon}$,
  where $\gamma,\delta,\epsilon=\pm1$, and a new rule is obtained:
 $x^{\gamma}z^{\delta}y^{\epsilon} \rightarrow  y^{\epsilon}x^{\gamma}z^{\delta}$.
 So, assumptions (A)  and (C) hold.

  We need to  check that the  assumption (B) holds also:\\
  Let $y^{\epsilon}z'^{\beta} \rightarrow  z'^{\beta}y^{\epsilon}$ and
 $z'^{\beta}x'^{\alpha} \rightarrow  x'^{\alpha}z'^{\beta}$,
  where $\alpha,\beta,\epsilon=\pm1$, be rules in $\Re_{0}$ whose overlap produced the  rule:
 $y^{\epsilon}x'^{\alpha}z'^{\beta} \rightarrow  z'^{\beta}y^{\epsilon}x'^{\alpha}$.
 We need to show that the overlap of $x^{\gamma}z^{\delta}y^{\epsilon} \rightarrow  y^{\epsilon}x^{\gamma}z^{\delta}$
 at left with $y^{\epsilon}x'^{\alpha}z'^{\beta} \rightarrow  z'^{\beta}y^{\epsilon}x'^{\alpha}$ at right
 resolve.
But,  there is an overlap of $x^{\gamma}z^{\delta}y^{\epsilon} \rightarrow  y^{\epsilon}x^{\gamma}z^{\delta}$
 at left with  $y^{\epsilon}z'^{\beta} \rightarrow  z'^{\beta}y^{\epsilon}$ at right which created the following rule: $x^{\gamma}z^{\delta}z'^{\beta}y^{\epsilon} \rightarrow  y^{\epsilon}x^{\gamma}z^{\delta}z'^{\beta}$.\\
  So, we have  that the overlap resolves (which is denoted by  $\checkmark$):\\
  $$\begin{array}{ccccccc}
   & \overbrace{x^{\gamma}z^{\delta}y^{\epsilon}}x'^{\alpha}z'^{\beta}\\
   \swarrow &&   \searrow & \\
 y^{\epsilon}x^{\gamma}z^{\delta}x'^{\alpha}z'^{\beta}&\checkmark&      \underbrace{x^{\gamma}z^{\delta}z'^{\beta}y^{\epsilon}}x'^{\alpha}
&  \rightarrow
& y^{\epsilon}x^{\gamma}z^{\delta}\underbrace{z'^{\beta}x'^{\alpha}}
 \end{array}$$
 So, at the first step the assumption (B) is satisfied.

Assume that at the $(n-1)$-th step of the Knuth-Bendix algorithm
of completion,  each new rule obtained satisfies the induction
assumptions  (A), (B) and (C). At the $n$-th step, a new rule is
obtained from an overlap ambiguity between a rule obtained at the
$(n-1)$-th step and a rule obtained at the $(n-1)$-th step or at
earlier steps.
Let $ux \rightarrow xu$ be a rule obtained at the $(n-1)$-th step which satisfies the condition
on  prefixes,  where
$u=t_{1}t_{2}..t_{k-1}t_{k}$ and $x,t_{i} \in
X\bigcup X^{-1}$.
Assume the rule $xvy \rightarrow yxv$  in $\Re_{n-1}$  satisfies the condition
on  prefixes, where
 $v$ is a word  and $y \in
X\bigcup X^{-1}$. \\
If $v$ is not the empty word, then by the induction assumption (B)
the overlap between  $ux \rightarrow xu$ and $xvy \rightarrow yxv$
resolve. So, let $v$ be the empty word, that is we consider the
overlap between the rules $ux \rightarrow xu$ and $xy \rightarrow
yx$:
$$\begin{array}{ccccccc}
 &  & uxy\\
 & \swarrow &  & \searrow & \\
xuy&  &  &  & uyx\\
\end{array}$$\\
There are three cases to check. We show first that in each case conditions (A) and (C) are satisfied.\\
 \textbf{First case:} If $xuy$ and $uyx$ are freely reduced and
 irreducible modulo $\Re_{n-1}$, then a new rule
is created: $uyx \rightarrow xuy$ which satisfies the condition on prefixes, that is the assumption (A) holds in
$\Re_{n}$ . Since by the induction assumption (C)  $\ell(ux)=\ell(xu)=(n-1)+2=n+1$, we have that
$\ell(uyx)=\ell(xuy)=(n+1)+1=n+2$, so (C) holds.\\
 \textbf{Second case:} Assume $xuy$ and $uyx$ are freely reduced but
 may be reducible modulo  $\Re_{n-1}$.\\
If there is a rule  $uy \rightarrow yu$ in $\Re_{n-1}$, then  the overlap resolves, by using first the rule $uy \rightarrow yu$ and then the rules $ux \rightarrow xu$ and $xy \rightarrow yx$.\\
If $u=u_{1}u_{2}$, where $u_{1}$ is a non-empty word,  and  $u_{2}y \rightarrow yu_{2}$ occurs in $\Re_{n-1}$, then we
have: $$\begin{array}{ccccccc}
 &  & \overbrace{ux}y=u_{1}u_{2}\overbrace{xy}\\
 &  \swarrow &&   \searrow & \\
&xu_{1}\underbrace{u_{2}y}&  &     u_{1}\underbrace{u_{2}y}x\\
&\downarrow && \downarrow\\
&xu_{1}yu_{2} &&u_{1}yu_{2}x \\
 \end{array}$$\\
 So,  a new rule is created:
 $u_{1}yu_{2}x \rightarrow xu_{1}yu_{2}$ which satisfies the condition on  prefixes, so
the conditions (A) and (C) are satisfied.

\textbf{Third case:} Assume $xuy$ and $uyx$ are not freely
reduced, then $uy$ is not freely reduced that is $t_{k}=y^{-1}$, since all the  subwords $xu$, $u$ and $yx$ are reduced.\\
 $$\begin{array}{ccccccc}
 &  & \overbrace{ux}y=t_{1}t_{2}..t_{k-1}y^{-1}\overbrace{xy}\\
 &  \swarrow &&   \searrow & \\
&xt_{1}t_{2}..t_{k-1}\underbrace{y^{-1} y}&  &     t_{1}t_{2}..t_{k-1}\underbrace{y^{-1}y}x\\
&\downarrow && \downarrow\\
&xt_{1}t_{2}..t_{k-1} &\checkmark &t_{1}t_{2}..t_{k-1}x\\
\end{array}$$\\
The word  $t_{1}t_{2}..t_{k-1}$ belongs to $\operatorname{pref}(u)$ and the
 rule $ux\rightarrow xu$ in $\Re_{n-1}$ satisfies the condition on prefixes, so by lemma
 \ref{lem:petitrule impliq}, there is a rule  $t_{1}t_{2}..t_{k-1}x \rightarrow
xt_{1}t_{2}..t_{k-1}$ which was obtained at an earlier step. So,
this overlap resolves using this rule.\\
We have shown that the conditions (A) and (C) hold in $\Re$.
Now, we need to show that the condition (B) holds too. \\
Let $ux
\rightarrow xu$ be a rule in $\Re_{n}$ which  satisfies the condition on  prefixes.
Assume $xvy \rightarrow yxv$ is  a
rule in $\Re_{n} \setminus \Re_{0}$, where $v=s_{1}s_{2}..s_{k}$ and $s_{i} \in X \bigcup X^{-1}$,
such that there is an overlap with $ux
\rightarrow xu$  at left and $xvy \rightarrow yxv$ at right.
In order to show that the condition (B) holds, we need to check that this overlap resolves.\\
From  assumption (A),  $xvy \rightarrow yxv$  satisfies also the condition on  prefixes.
Then:

$$\begin{array}{ccccccc}
 &  & \overbrace{ux}vy=u\overbrace{xs_{1}s_{2}..s_{k}y}\\
 &  \swarrow &&   \searrow & \\
&xuvy=xus_{1}s_{2}..s_{k}y&  &    uyxv=uyxs_{1}s_{2}..s_{k}
\end{array}$$

Since $xvy \rightarrow yxv$  satisfies the condition on  prefixes,  there are rules $xy \rightarrow yx$ and
$ys_{i}\rightarrow s_{i}y$ for $1\leq i\leq k$  in $\Re_{0}$. So,
there is an overlap of the rule $ux \rightarrow xu$ at left and
$xy \rightarrow yx$ at right, which created the rule $uyx
\rightarrow xuy$.
So, we have using this rule and the rules $ys_{i}\rightarrow
s_{i}y$,  for $1\leq i\leq k$ that:\\
 $$\begin{array}{ccccccc}
  \overbrace{ux}vy=u\overbrace{xs_{1}s_{2}..s_{k}y}\\
   \swarrow    &\searrow \\
xuvy & \checkmark    \underbrace{uyx}s_{1}..s_{k}
 \rightarrow
xu\underbrace{ys_{1}}..s_{k}
 \rightarrow
  ..\rightarrow
 xus_{1}..s_{k}y
\end{array}$$\\
So, this overlap resolves.
At last, we need to check that the overlaps of rules with rules of free reduction resolve.
Let  $ux \rightarrow xu$ be a rule  which satisfies the condition
on prefixes,  where
$u=t_{1}t_{2}..t_{k-1}t_{k}$ and $x,t_{i} \in
X\bigcup X^{-1}$. We check first the overlap with $ux \rightarrow xu$ at left and  $xx^{-1}\rightarrow 1$ at right:
$$\begin{array}{ccccccc}
  \overbrace{ux}x^{-1}=u\overbrace{xx^{-1}} \\
   \swarrow    & \searrow \\
x\underbrace{ux^{-1}}  \rightarrow
\underbrace{xx^{-1}}u
 \rightarrow u &\checkmark& u
\end{array}$$\\

We check now the overlap with $t^{-1}_{1}t_{1}\rightarrow 1$  at left and $ux \rightarrow xu$ at right.
$$\begin{array}{ccccccc}
  \overbrace{t_{1}^{-1}t_{1}}t_{2}..t_{k-1}t_{k}x=t_{1}^{-1}\overbrace{ t_{1}t_{2}..t_{k-1}t_{k}x}\\
   \swarrow    &\searrow \\
t_{2}..t_{k-1}t_{k}x &\checkmark     \underbrace{t_{1}^{-1}x}u
 \rightarrow
x\underbrace{t_{1}^{-1}u}
 \rightarrow
\underbrace{xt_{2}}..t_{k-1}t_{k}
 \rightarrow
  ..
\end{array}$$\\
(ii) From (i), all the rules in $\Re$ have the form $ux \rightarrow xu$.
So, the word $ux$ is positive if and only if the word $xu$ is positive.
\end{proof}

\begin{prop}\label{prop:RsatisfiesR+}
Let  $\Re$ be the  complete rewriting system obtained from the
Knuth-Bendix algorithm  of completion of $\Re_{0}$.
Then  $\Re^{+}$ is a rewriting system for the right-angled Artin
monoid $M$ defined by $\Im $.
\end{prop}
\begin{proof}
We will show by induction on the number of steps in the
Knuth-Bendix algorithm
 of completion of $\Re_{0}$ that  at each step a positive rule is
 created by the overlap of two positive rules.
At the $0-$th step, the set of positive rules is  $\Re_{0}^{+}=\{
x^{\gamma}y^{\epsilon} \rightarrow y^{\epsilon}x^{\gamma}$,  where
$\gamma=\epsilon=1$, for each edge $ x-y \in \Im \}$.\\
The rewriting system $\Re_{0}^{+}$ is a  rewriting system for
the right-angled Artin monoid $M$ defined by $\Im $, since to each edge in $\Im $, there corresponds exactly one rule in  $\Re_{0}^{+}$. \\
At the first step, we have overlaps between rules
 $x^{\gamma}y^{\epsilon} \rightarrow  y^{\epsilon}x^{\gamma}$ and
 $y^{\epsilon}z^{\delta} \rightarrow  z^{\delta}y^{\epsilon}$,
   where $\gamma,\delta,\epsilon=\pm1$, and the new rule obtained is
 $x^{\gamma}z^{\delta}y^{\epsilon} \rightarrow  y^{\epsilon}x^{\gamma}z^{\delta}$.
The  rule  $x^{\gamma}z^{\delta}y^{\epsilon} \rightarrow
y^{\epsilon}x^{\gamma}z^{\delta}$ is a positive rule if and only
$\gamma=\delta=\epsilon=1$. Note that the words
$x^{\gamma}z^{\delta}y^{\epsilon}$ and
$y^{\epsilon}x^{\gamma}z^{\delta}$ are freely reduced, since
otherwise $z=x$ and this would contradict the total ordering on
the generators.\\
Assume that at the $(n-1)-$th step  in  the Knuth-Bendix algorithm
of completion, all the positive rules were created by  overlaps
between positive rules. We will assume that at  the $n-$th step,
there is a positive rule which is  created by an overlap between
rules which are not positive. \\
Let $ux\rightarrow xu$ be a rule in $\Re_{n-1}$, which satisfies the condition on  prefixes,  where
$u=t_{1}t_{2}..t_{n-1}t_{n}$ and
 $x, t_{i} \in X\bigcup X^{-1}$.\\
  Let  $xy \rightarrow yx$ be a rule in
 $\Re_{0}$, where $y \in X\bigcup X^{-1}$, such that there is an overlap of $ux\rightarrow xu$
 at left and $xy \rightarrow yx$ at right which creates a positive rule.\\
By proposition \ref{prop:Rsatisfait_pos}(ii),  $\Re$ satisfies the following condition:  one side
 of a rule in $\Re$ is a positive word  if and only if the other side is a positive word.
 So, if $ux$ or $xy$ are not positive words, then the words $uyx$ and $xuy$ are not positive also and
  $uyx$ and $xuy$ reduce to  positive words only if there  are
  free reductions which eliminate the generators with negative
  exponent sign.
  So,  $t_{n}=y^{-1}$ and we have:
 $$\begin{array}{ccccccc}
   & \overbrace{ux}y=t_{1}t_{2}..t_{n-1}y^{-1}\overbrace{xy}\\
   \swarrow &&  \searrow & \\
xuy=xt_{1}t_{2}..t_{n-1}\underbrace{y^{-1}y}&&      t_{1}t_{2}..t_{n-1}\underbrace{y^{-1}y}x\\
\downarrow && \downarrow\\
xt_{1}t_{2}..t_{n-1} && t_{1}t_{2}..t_{n-1}x\\
\end{array}$$\\
The word  $t_{1}t_{2}..t_{n-1}$ belongs to $\operatorname{pref}(u)$ and the
 rule $ux\rightarrow xu$  satisfies the condition on  prefixes, so from lemma
 \ref{lem:petitrule impliq}, there is a rule  $t_{1}t_{2}..t_{n-1}x \rightarrow
xt_{1}t_{2}..t_{n-1}$ and from the condition (C) in proposition
\ref{prop:Rsatisfait_pos} it was obtained at the $(n-2)-$th step.
So, a positive rule cannot be created by this overlap.
\end{proof}
\textbf{Proof of theorem \ref{thm:artin_embed}}
\begin{proof}
Let define a total ordering $<$ on the generators and their inverses and
let $\Re_{0} =\{ x^{\gamma}y^{\epsilon} \rightarrow
y^{\epsilon}x^{\gamma}$,  where
$\gamma$ and $\epsilon$ take all values $\pm 1$, for each edge $ x-y \in \Im \}$ be a  rewriting system for the right angled Artin group defined by $\Im$, using the length-lexicographical ordering induced by $<$. Let $\Re$ be the complete rewriting system  obtained from the Knuth-Bendix algorithm of completion applied  on $\Re_{0}$. \\
Then by proposition \ref{prop:Rsatisfait_pos}, $\Re$ satisfies the condition $C^{+}$  and from proposition \ref{prop:RsatisfiesR+}, $\Re$ satisfies the condition on $\Re^{+}$.\\
 So, from theorem \ref{thm_result}, the right angled Artin monoid embeds in the corresponding group.
\end{proof}
 In  appendix B, we  give a very simple proof of
the embedding of right angled Artin monoids  with  complete
defining graph or with defining graph, such that all the closed
paths have even length, in their corresponding groups. In fact, we
find in these two cases a finite and complete rewriting system
(with no need of completion) which satisfies the condition $C^{+}$
and the condition on $\Re^{+}$, by using some elementary tools
from graph theory.
\section{Appendix A
:proof of lemma \ref{lem:petitrule impliq}}
 Before we give the proof, we need the
following lemma and definition.
\begin{lem}
Assume the set $\{t_{1},t_{2},..,t_{k}\}$ is a clique. Then any
word $w=t_{i_{1}}t_{i_{2}}..t_{i_{k}}$, where $i_{j} \in
\{1,2,..,k\}$, can be reduced to its normal form modulo $\Re$ by
using only rules from $\Re_{0}$. That is, no new rule is created
from the overlaps between rules involving only elements from a
clique.
\end{lem}
\begin{proof}
Since there is a total ordering on the generators, each word
$w=t_{i_{1}}t_{i_{2}}..t_{i_{k}}$, where $i_{j} \in \{1,2,..,k\}$
reduces to a unique irreducible word, which is the least one
lexicographically and which is obtained by the application of rules
from $\Re_{0}$.
\end{proof}
Let $w$ be a word in the free group generated by $X$. Assume
$w=t_{1}t_{2}t_{3}...t_{k}$, where $t_{i} \in X\bigcup X^{-1}$ for
$1\leq i \leq k$. Then, we say that $t_{1},t_{2},..,t_{i-1}$
\textbf{ appear before} $t_{i}$ and
$t_{i+1},t_{i+2},..,t_{k}$  \textbf{appear after} $t_{i}$.
\begin{proof}
 We prove  that by induction on the
number of steps in the algorithm of completion
and we prove also the following condition (**):
if $x$ is the last letter in the left-hand side of a rule and $x$ commutes with all the letters appearing before, then in the right-hand side $x$ is the first letter.\\
At the $0-$th step, the assumptions hold trivially.
At the first step
of completion,  there is an overlap between
$t_{1}x \rightarrow xt_{1}$ and $xt_{i_{2}} \rightarrow
t_{i_{2}}x$ in $\Re_{0}$,  $i_{2}\neq 1$, which gives
$$\begin{array}{ccccccc}
 &  & t_{1}x t_{i_{2}}\\
 & \swarrow &  & \searrow & \\
xt_{1} t_{i_{2}}&  &  &  & t_{1} t_{i_{2}}x\\
\end{array}$$\\
The rule $t_{1} t_{i_{2}}x \rightarrow xt_{1} t_{i_{2}}$ is obtained and assumption (**) holds.
The word  $xt_{1} t_{i_{2}}$ is irreducible modulo $\Re$ since it
is a subword of  $xt_{1}t_{i_{2}}..t_{i_{n}}$. The word  $t_{1} t_{i_{2}}x$ is irreducible modulo
$\Re_{0}$,  since each subword is irreducible modulo $\Re_{0}$:
the word $t_{1} t_{i_{2}}$ is irreducible modulo $\Re_{0}$ since
otherwise there would be an inclusion ambiguity with the rule
$t_{1}t_{i_{2}}..t_{i_{n}}x \rightarrow
xt_{1}t_{i_{2}}..t_{i_{n}}$ in  $\Re$ and this is not possible.\\
Assume that at the $(m-1)-th$  step of
completion of $\Re_{0}$, the rule $t_{1}t_{i_{2}}..t_{i_{n-2}}x
\rightarrow xt_{1}t_{i_{2}}..t_{i_{n-2}}$ is obtained, where
$t_{i_{j}}\in \{t_{2},..,t_{k}\}$. Then at the $m-th$  step, from the overlap between
this rule and $xt_{i_{n-1}} \rightarrow t_{i_{n-1}}x$ ($i_{n-1}\neq 1$), we have: $$\begin{array}{ccccccc}
 &  & \overbrace{t_{1}t_{i_{2}}..t_{i_{n-2}}x} t_{i_{n-1}}\\
 & \swarrow &  & \searrow & \\
xt_{1}t_{i_{2}}..t_{i_{n-2}}t_{i_{n-1}}&  &  &  & t_{1}t_{i_{2}}..t_{i_{n-2}}t_{i_{n-1}}x\\
\end{array}$$\\
We need to show that a new rule is created, with the word $xt_{1}t_{i_{2}}..t_{i_{n-2}}t_{i_{n-1}}$ in the right-hand side and the word $t_{1}t_{i_{2}}..t_{i_{n-2}}t_{i_{n-1}}x$ in the left-hand side. To do that, we need to show that both these words are irreducible modulo the rewriting system $\Re_{m-1}$ obtained one step before.
In fact, the word $xt_{1}t_{i_{2}}..t_{i_{n-2}}t_{i_{n-1}}$ is irreducible
modulo $\Re$ since it is a subword of
$xt_{1}t_{i_{2}}..t_{i_{n}}$, which is an irreducible
word modulo $\Re$. We show that the word
$t_{1}t_{i_{2}}..t_{i_{n-2}}t_{i_{n-1}}x$ is  irreducible modulo
$\Re_{m-1}$, by showing that each of its subword cannot occur as a left-hand side of a rule in $\Re_{m-1}$.
If there is a
rule $t_{i_{j}}t_{i_{j+1}}..t_{i_{n-1}}x\rightarrow ..$, where
$i_{j} \neq 1$. Then from assumption (**), we
have that $t_{i_{j}}t_{i_{j+1}}..t_{i_{n-1}}x\rightarrow x..$,
since $x$ commutes with all the $t_{i_{j}}$ and this would imply
that $t_{i_{j}}>x$. But from the existence of the rule $xt_{i_{j}}
\rightarrow t_{i_{j}}x$ in $\Re_{0}$, it holds that $x>t_{i_{j}}$.
So, there can be no such rule and it is easy to check that for the other subwords, so the word
$t_{1}t_{i_{2}}..t_{i_{n-2}}t_{i_{n-1}}x$ is  irreducible modulo
$\Re_{m-1}$.\\
So, there is the   rule $t_{1}t_{i_{2}}..t_{i_{n-2}}t_{i_{n-1}}x
\rightarrow xt_{1}t_{i_{2}}..t_{i_{n-2}}t_{i_{n-1}}$ in $\Re$ and assumption (**) holds.\\
From the definition of $\Re_{0}$, if there are rules  $t_{1}x \rightarrow xt_{1}$ and
  $xt_{i} \rightarrow t_{i}x$ for
 $2 \leq i \leq k$, then there are also the rules  $t_{1}x^{-1} \rightarrow x^{-1}t_{1}$ and
  $x^{-1}t_{i} \rightarrow t_{i}x^{-1}$ for
 $2 \leq i \leq k$. So, in the same process as before the rule $ux^{-1} \rightarrow x^{-1}u$
 is created.
\end{proof}

\section{Appendix B}
\subsection{Right angled Artin groups with  defining graph such that
all the closed paths have even length} Let $\Im$ be the defining
graph of a right-angled Artin group $G$  and the corresponding right-angled Artin monoid $M$,  such that all the closed
paths in $\Im$ have even length. As before, we denote by $X$ the set of
vertices in $\Im$ or equivalently the set of generators of $G$ and
by $X^{-1}$ the set of inverses of the elements in $X$.
 We will show that $G$ admits a finite and complete rewriting system $\Re$  satisfying
  the condition $C^{+}$ by using some very simple tools from graph theory (see \cite{wilson}).
In fact, we will define  a total ordering $'>'$ on $X \bigcup
X^{-1}$  such that the rewriting system \\
 $\Re =\{x_{i}^{\gamma}x_{j}^{\epsilon} \rightarrow
x_{j}^{\epsilon}x_{i}^{\gamma}$,  where $x_{i}'>' x_{j}$ and
$\gamma$, $\epsilon$ take all values $\pm 1\}$ is finite and
complete with no need of completion.
\begin{lem}
Let $\Im$ be the defining graph of a right-angled Artin group such
 that  all the closed paths  have even length.
 Then the vertices of $\Im$ can be coloured in two colours (black and white) such that any two
  adjacent vertices (i.e vertices
 connected by an edge) have different colours.
\end{lem}
\begin{proof}
The graph $\Im$ is simple and satisfies the condition that all the
closed paths have even length, so using the proposition
\cite[p.81]{wilson} the vertices of  $\Im$ can be coloured in two
colours: black and white.
\end{proof}
Assume a colouring is assigned to $\Im$. Denote by $B$ the set of
generators corresponding to vertices  which are coloured in black
in $\Im$ and by $W$ the set of generators corresponding to
vertices  which are coloured in white  in $\Im$. By definition of the colouring,  $B\bigcap W=\emptyset$ and
$X=B \bigcup W$. We denote by $B^{-1}$ the set of inverses of the
elements in $B$  and by $W^{-1}$
 the set of inverses of the elements in $W$.
\begin{lem}\label{lem:tot_order}
The colouring of the vertices in $\Im$ defines a well-founded,
total  ordering on $X \bigcup X^{-1}$.
\end{lem}
\begin{proof}
We order  arbitrarily the elements in $B\bigcup B^{-1}$ and also
the elements in  $W\bigcup W^{-1}$, but we require  that the least
element from  $B\bigcup B^{-1}$ is greater than the greatest
element from  $W\bigcup W^{-1}$.\ This ordering is well-founded,
i.e there can be no infinite sequence \newline
$x_{1}>x_{2}>..>x_{k}>...$, since $(B \bigcup B^{-1}) \bigcap (W
\bigcup W^{-1})= \emptyset$ and it is  total  since  $X\bigcup
X^{-1}=(B \bigcup W )\bigcup (B^{-1} \bigcup W^{-1})$.
\end{proof}
\begin{rem}
From the proof of lemma\ref{lem:tot_order}, it results for each
edge $a - b$ that the corresponding rules are  oriented
$a^{\gamma}b^{\epsilon}\rightarrow b^{\epsilon}a^{\gamma}$ (where
$\gamma$ and $\epsilon$ take all values $\pm1$) if and only if
 $a$ is coloured in black and $b$ in white.
\end{rem}
\begin{prop}
Let $>_{lex}$ be the length-lexicographical ordering induced by
the total ordering on $X \bigcup X^{-1}$ defined above and let
$\Re$ be the rewriting system obtained using $>_{lex}$. Then $\Re$
is finite and complete and  it satisfies  the condition $C^{+}$.
Further,   $\Re^{+}$ is a rewriting system for $M$.
\end{prop}
\begin{proof}
$\Re$ is terminating since it is defined by a
length-lexicographical ordering, it is finite  and it satisfies
$C^{+}$ from lemma \ref{lem:order_resolve}. So, it remains to show
that during the Knuth-Bendix algorithm of completion, the overlaps
of the rules resolve with no need to add new rules. In fact, we
will show that there are no overlaps between the rules themselves
but only between the rules and the rules  of kind
 $xx^{-1} \rightarrow 1$ and $x^{-1}x \rightarrow1$ ($x \in X$) and these resolve from lemma \ref{lem:order_resolve}.
 Assume there is an overlap between two rules $x^{\gamma}y^{\epsilon}\rightarrow y^{\epsilon}x^{\gamma}$ and $y^{\epsilon}z^{\delta}\rightarrow z^{\delta}y^{\epsilon}$,
  where $x,y,z \in X$ and $\gamma, \delta, \epsilon =\pm 1$. If there is a rule
   $x^{\gamma}y^{\epsilon}\rightarrow y^{\epsilon}x^{\gamma}$ in $\Re$, then it means that
  there is an edge $x-y$ in $\Im$ with $x$ coloured in black and $y$ coloured in white. In the same way,
   if there is a rule $y^{\epsilon}z^{\delta}\rightarrow z^{\delta}y^{\epsilon}$ in $\Re$, then it means that
  there is an edge $y-z$ in $\Im$ with $y$ coloured in black and $z$ coloured in white. But this is a contradiction to the colouring of $\Im$.
\end{proof}
\subsection{Right angled Artin groups with a complete defining graph}
Let $M$ be the free commutative monoid on $n$ generators with the
standard presentation  $M=\langle x_{1},..,x_{n}\mid
x_{i}x_{j}=x_{j}x_{i},n\geq i>j\geq 1>$. The left and right graphs
$\Gamma_{l}$ and $\Gamma_{r}$ both have cycles, so nothing can be
concluded here using Adian's
criteria.\\
Yet, it is well known that this monoid  embeds in the free abelian
group on $n$ generators or in other words in the right angled
Artin group with a complete defining graph on $n$ vertices.
 We  show this well-known result, using the criteria of section 3. Indeed,
we  give a complete and finite rewriting system which
satisfies the condition $C^{+}$.\\
Let $\langle x_{1},..,x_{n},x_{1}^{-1},..,x_{n}^{-1}\mid
x_{i}x_{j}=x_{j}x_{i},n\geq i>j\geq 1>,
x_{i}^{-1}x_{i}=x_{i}x_{i}^{-1}=1, 1\leq i \leq n\rangle  $ be the monoid
presentation of the free abelian group on $n$ generators. Let
$\Re$ be the rewriting system defined for the presentation above,
using a length-lexicographical ordering induced by the following
total ordering on the generators:\
 $x _{n}> x_{n}^{-1}>...>x_{2} > x_{2}^{-1}> x_{1} > x_{1}^{-1}$.\
  So,  $\Re$  can be described as the following set of rules:
$\Re =\{ x_{i}^{\gamma}x_{j}^{\epsilon} \rightarrow
x_{j}^{\epsilon}x_{i}^{\gamma}$,  where $n\geq i>j\geq 1$ and
$\gamma, \epsilon=\pm 1 \}$
This rewriting system is terminating since it is defined by a
length-lexicographical ordering and one can check that the
overlaps between the rules which occur resolve with no need to add
new relations. So $\Re$ is finite and complete and it satisfies
the condition $C^{+}$ and the condition on  $\Re^{+}
$. By theorem \ref{thm_result}, the monoid
embeds in the group.

\end{document}